\documentclass[12pt,a4paper]{amsart}
\usepackage{a4}
\usepackage[english]{babel} 
\usepackage{amsmath,amssymb,amsthm}
\def\leq {\leqslant}
\def\le {\leqslant}

\def\geq {\geqslant}
\def\@bibitem[#1]#2{\item\@biblabel{#1}.\if@filesw
{\def\protect##1{\string##1\space}\immediate\write
\@auxout{\string\bibcite{#2}{#1}}}\fi\ignorespaces\@showtag{#2}}

\theoremstyle{plain}
\newtheorem{theorem}{Theorem}
\newtheorem{rem}{Remark}
\newtheorem{lemma}{Lemma}
\renewcommand{\theequation}%
{\arabic{section}.\arabic{equation}}
\begin{document}

\title{ On estimates of the order of approximation of functions of several variables in the anisotropic Lorentz - Karamata space}
\author{ G. Akishev}
\address{ Lomonosov Moscow University, Kazakhstan Branch \\
Str. Kazhymukan, 11 \\
010010, Nur-Sultan, Kazakhstan}

\address{
Institute of mathematics and mathematical modeling\\
Pushkin str, 125 \\
050010, Almaty \\
 Kazakhstan
 }
 
\maketitle

\begin{quote}
\noindent{\bf Abstract.}
In this paper we consider  anisotropic Lorentz-Karamata space $2\pi$ of periodic functions of $m$ variables and Nikol'skii--Besov's class .  
In this paper, we establish order-sharp estimates of the best approximation by trigonometric polynomials with harmonic numbers from the step hyperbolic cross of functions from the Nikol'skii - Besov class in the norm of the anisotropic Lorentz-Karamata space.
\end{quote}
\vspace*{0.2 cm}



\section*{Introduction}

Let $\mathbb{N}$, $\mathbb{Z}$, $\mathbb{R}$ be sets of natural, integer, real numbers, respectively, and $\mathbb{Z}_{+} = \mathbb{N}\cup\{0\}$,
$\mathbb{R}^{m}$ --- $m$ -- dimensional Euclidean point space $\bar{x} = (x_{1}, \ldots, x_{m})$ with real coordinates; $I^{m} = \{\bar{x} \in \mathbb{R}^{m}; \ 0 \leq x_{j} \leq 1; \ j = 1, \ldots, m \}=[0, 1]^{m} $ --- $m$--dimensional
cub and $\mathbb{Z}_{+}^{m}$ is the Cartesian product of the sets $\mathbb{Z}_{+}$, $m$ is times.

 \smallskip

\mbox{D\ e\ f\ i\ n\ i\ t\ i\ o\ n} (see \cite[p. 10, 51]{1}, \cite[p. 6 , 24]{2}, \cite[p. 108]{3}).  
\ \ A positive and Lebesgue measurable function $v(t)$ is called slowly varying by $[1, + \infty)$ if, for any $\varepsilon> 0$, the function $t^{\varepsilon} v(t)$ is equivalent to a function non-decreasing on $[1, \infty)$
   and the function $t^{-\varepsilon} v(t)$ is equivalent to a function non-increasing on $[1, \infty)$.

The set of such functions is denoted by $SV[1, \infty)$. 
    \smallskip

Examples of slowly varying functions (see \cite{1}, \cite[p. 108]{3}):

 1. $ l_1(t)=1+\log t \in SV [1,\infty ), \ \ 
 l_2 (t)=1+\log l_1(t),\  l_3(t)=1+\log l_2(t).$

  2. $l_i(t)=1+l_{i-1}(t)\in SV [1,\infty), i=2,3,\ldots s.$

Here and below the symbol $\log t = \log_{2} t$ for $t> 0$.

For a given slowly varying function $v$ on $ [1, \infty)$, we put $V(t) = v(1/t)$ for $t \in (0, 1]$.

Let numbers $p \in (1, \infty)$, $\tau\in [1, \infty)$ and a function $v\in SV[1, \infty)$. Lorentz-Karamata space  $L_{p, V,  \tau}(\mathbb{T})$ is the set of all Lebesgue measurable and $2 \pi$ periodic functions $f$ for which (see \cite[p. 112]{3})
$$
\|f\|_{p, V, \tau}^{*} := \left\{\int_{0}^{1}\bigl(f^{*}
(t)\bigr)^{\tau}V^{\tau}(t)t^{\frac{\tau}{p}-1}dt
\right\}^{\frac{1}{\tau}} < +\infty,
$$
where $f^{*}(t)$ is a non-increasing rearrangement of the function
$|f(2\pi x)|$, $x \in [0, 1]$, $\mathbb{T}=[0, 2\pi]$ (see e.g. \cite{3}) .

Note that for $V(t) = 1$ the space $L_{p, V, \tau }(\mathbb{T})$ coincides with the Lorentz space $L_{p, \tau}(\mathbb{T})$, which consists of all functions $f$ such that (see \cite[p. 228]{4})
$$
\|f\|_{p,\tau}^{*}:=\Bigg(\frac{\tau}{p}\int\limits_{0}^{1}\bigl(f^{*}(t)\bigr)^{\tau}t^{\frac{\tau}{p}-1}dt
\Bigg)^{1/\tau} < \infty,
$$
for $1 < p < \infty,$ $1\leqslant\tau<+\infty$.

It is known that the Lorentz - Karamata space $L_{p, V, \tau}(\mathbb{T})$
is one example of a symmetrical space (see \cite[p. 113]{3}).

Let 
 $\bar{p}=(p_{1},\ldots p_{m})$, $\bar{\tau}=(\tau_{1},\ldots \tau_{m}),$  and $p_{j} \in (1, \infty)$, $\tau_{j}\in [1, \infty)$ and functions $v_{j}\in SV[1, \infty)$, $j=1,...,m$.
By $L_{\overline{p}, \overline{V}, \overline{\tau}}^{*}(\mathbb{T}^{m})$ we denote the anisotropic Lorentz-Karamata space of all measurable Lebesgue functions of $m$ variables $f$ with period
$2 \pi$ for each variable and for which the quantity (see \cite{5}) 
$$
\|f\|_{\bar{p}, \bar{V}, \bar{\tau}}^{*} := \|\ldots\|f^{*_{1},...,*_{m}}\|_{p_{1}, V_{1}, \tau_{1}}\ldots\|_{p_{m}, V_{m}, \tau_{m}}
$$
$$
\Bigl[\int_{0}^{1}
\Bigl[\ldots\Bigl[\int_{0}^{1}\left(
f^{*_{1},...,*_{m}}(t_{1},...,t_{m})
\right)^{\tau_{1}}\biggl(\prod_{j=1}^{m}V_{j}(t_{j})t_{j}^{\frac{1}{p_{j}}-\frac{1}{\tau_{j}}}\biggr)^{\tau_{1}}dt_{1}\Bigr]^
{\frac{\tau_{2}}{\tau_{1}}} \ldots
\Bigr]^{\frac{\tau_{m}}{\tau_{m-1}}}dt_{m}
\Bigr]^{\frac{1}{\tau_{m}}}
$$
is finite, where
$f^{*_{1},...,*_{m}}(t_{1},...,t_{m})$
non-increasing rearrangement  of a function
$|f(2\pi \bar{x})|$ for each variable $x_{j} \in [0, 1]$ with fixed other variables (see \cite{5}, \cite{6}). Here and in what follows, $\bar{V}=\bar{V}(t)=(V_{1}(t),...,V_{m}(t))$ and  $\mathbb{T}^{m} = [0, 2\pi]^{m}$.
 
As usual, the space $l_{\overline{p}}$ consists of sequences of real numbers $ \left\{a_{\overline{n}}\right\}_{\overline{n} \in \mathbb{Z}_{+}^{m}}$ such that
\begin{equation*}  
\Bigl\|
\Bigl\{a_{\overline{n}}\Bigr\}_{\bar{n} \in \mathbb{Z}_{+}^{m}}
\Bigr\|_{l_{\overline{p}}(\mathbb{Z}_{+}^{m})} = \Bigl\{\sum\limits_{n_{m}
=0}^{\infty }\Bigl[ ...\Bigl[\sum\limits_{n_{1} =0
}^{\infty }\Bigl|a_{\overline{n}}\Bigr|^{p_{1}}\Bigr]
^{\frac{p_{2}}{p_{1}}} ...\Bigr]^{\frac{p_{m}}{p_{m-1}}}
\Bigr\}^{\frac{1}{p_{m}}} <+\infty ,
\end{equation*} 
where $\overline{p} =\left( p_{1} ,...,p_{m} \right)$, $1\leq
p_{j} <+\infty$, $j=1,2,...,m$.

If $p_{j}=\infty, \;\; j=1,...,m$, then  
$$
\Bigl\|\{a_{\bar{n}}\}\Bigr\|_{l_{\infty}(\mathbb{Z}_{+}^{m})}=\sup\limits_{\bar{n}\in\mathbb{Z}_{+}^{m}}|a_{\bar{n}}|.
$$

We will use the following notation:
let $\overset{\circ \;\;}
L_{\bar{p}, \bar{V}, \bar{\tau}}^{*}
(\mathbb{T}^{m})$ be the set of functions $f\in
L_{\bar{p}, \bar{V}, \bar{\tau}}^{*}(\mathbb{T}^{m})$ such
that
 $$
\int\limits_{0}^{2\pi }f\left(\overline{x} \right) dx_{j}
=0, \;\;   j=1,...,m;
 $$
$a_{\overline{n} } (f)$ be the Fourier coefficients  of the function $f\in
L_{1}(\mathbb{T}^{m})$ with respect to the  system $\{e^{i\langle\overline{n} ,\overline{x}\rangle}\}_{\overline{n}\in \mathbb{Z}^{m}}$ and $\langle\bar{y}, \bar{x}\rangle=\sum\limits_{j=1}^{m}y_{j}x_{j}$, 
$$
\delta_{\overline{s}}(f, \overline{x})
:=\sum\limits_{\overline{n} \in \rho(\overline{s})
}a_{\overline{n}}(f) e^{i\langle\overline{n} ,\overline{x}\rangle } ,
$$
 $$
\rho(\bar{s}):=\left\{ \overline{k} =( k_{1}
,...,k_{m}) \in \mathbb{Z}^{m}: \, \,  [2^{s_{j} -1}] \leq \left| k_{j}
\right| <2^{s_{j} } ,j=1,...,m\right\}, s_{j}\in \mathbb{Z}_{+},
$$ 
where $[a]$ is the integer part of the number $a \in \mathbb{Z}_{+}$.

Let $\overline{\gamma}=(\gamma_{1},\ldots,\gamma_{m}),$ $\gamma_{j}>0, j=1,\ldots,m$.
$E_{n}^{(\overline{\gamma})}(f)_{\overline{p}, \overline{V}, \overline{\tau}}$ is the best approximation of the function $f\in L_{\overline{p}, \overline{V}, \overline{\tau}}^{*}(\mathbb{T}^{m})$ trigonometric polynomials with harmonic numbers from the set $Q_{n}^{\overline\gamma}=
 \cup_{{}_{\langle\overline{s}, \overline{\gamma}\rangle < n}}\rho(\overline{s})$ is a stepped hyperbolic cross. For $\overline{V}=(1, \ldots , 1)$ and $\overline{\tau}=\overline{p}$ the  quantity  $E_{n}^{(\overline{\gamma})}(f)_{\overline{p}, \overline{V}, \bar{\tau}}$ is denoted as $E_{n}^{(\overline{\gamma})}(f)_{\overline{p}}$ (see \cite{12} and references therein).

We put  
In this article, we will consider an analogue of the Nikol'skii-Besov class in the anisotropic Lorentz-Karamata space:
 $$
S_{\overline{p}, \overline{V}, \overline{\tau}}^{\bar r}B=
\Bigl\{f\in \overset{\circ \;\;}{L}_{\overline{p}, \overline{V}, \overline{\tau}}^{*}(\mathbb{T}^{m}) : \quad
\|f\|_{\overline{p}, \overline{V}, \overline{\tau}}^{*} + \Bigl\|\Bigl\{\prod_{j=1}^{m}
2^{s_{j}r_{j}} \|\delta_{\bar
s}(f)\|_{\overline{p}, \overline{V}, \overline{\tau}}^{*} \Bigr\}_{\bar{s}\in
\mathbb{Z}_{+}^{m}}\Bigr\|_{l_{\bar\theta}}\leq 1\Bigr\},
 $$
where 
 $\bar{\theta}=(\theta_{1},...,\theta_{m}),$  $\bar{r}=(r_{1},...,r_{m}),$ 
$1\leq\theta_{j}\le+\infty,$ $0<r_{j}<+\infty,$  $j=1,...,m.$

In case $V_{j}(t) = 1 $ and $\tau_{j} = p_{j}$, $j = 1, \ldots, m $ , class $S_{\overline{p}, \overline{V}, \overline{\tau}}^{\overline{r}} B $ coincides with the well-known Nikol'skii-Besov class in the Lebesgue space (see  \cite{7} - \cite{9}.

Order-exact estimates $E_{n}^{(\overline{\gamma})}(f)_{\bar{p}}$ - the best approximation of functions from the Sobolev classes $W_{p}^{\bar r}$ and Nikol'skii - Besov $S_{p, \theta}^{\bar r}B$, in the metric of the space $L_{q}(\mathbb{T}^{m})$, $1< p, q < \infty$ are well known and are given in the survey articles \cite{10}, \cite{11} and in the monographs \cite{12} (for more details, see the references therein).
These questions in the anisotropic Lorentz space $L_{\bar{q}, \bar{\tau}}^{*}(\mathbb{T}^{m})$ are investigated in \cite{13}-\cite{18} and in the anisotropic generalized Lorentz space are studied in \cite{19}, \cite{20}.

We put 
$$
E_{n}^{(\overline{\gamma})}(S_{\bar{p}, \bar{V}^{(1)}, \bar{\tau}^{(1)},  \bar{\theta}}^{\bar r}B)_{\bar{q}, \bar{V}^{(2)}, \bar{\tau}^{(2)}} := \sup\limits_{f \in S_{\bar{p}, \bar{V}^{(1)}, \bar{\tau}^{(1)}, \bar{\theta}}^{\bar r}B}E_{n}^{(\overline{\gamma})}(f)_{\bar{q}, \bar{V}^{(2)}, \bar{\tau}^{(2)}},
$$
where $\bar{q}=(q_{1},\ldots , q_{m})$, $\bar{V}^{(i)}(t)=(V_{1}^{(i)}(t),\ldots , V_{m}^{(i)}(t))$, $t\in (0, 1]$, $\bar{\tau}^{(i)}=(\tau_{1}^{(i)},\ldots , \tau_{m}^{(i)})$ and $1<q_{j}, \tau_{j}^{(i)}< \infty$, $i=1, 2$,  $j=1,...,m$.

The main aim of the present paper is to find the order of the quantity
$$
E_{n}^{(\overline{\gamma})}(S_{\bar{p}, \bar{V}^{(1)}, \bar{\tau}^{(1)},  \bar{\theta}}^{\bar r}B)_{\bar{q}, \bar{V}^{(2)}, \bar{\tau}^{(2)}}.
$$

This paper is organized as follows. In  Section 1 we give auxiliary results.  
 In Section 2, we present and prove the main result.
 
  We shall denote by  $C(p,q,y,..)$  positive quantities which depend
only on the parameter in the parentheses and not necessarily the
same in distinct formulas .
 The notation $A\left( y\right) \asymp
B\left( y\right)$ means that there exist positive constants
 $C_{1},\,C_{2} $ such that  $C_{1} A(y)
\leq B(y) \leq C_{2} A(y)$.
 For brevity, in the case of the inequalities $ B \geq C_{1} A$ or $B \leq C_{2}A$, we often write $B \gg A$ or $B \ll A $, respectively.

\section{Auxilary statements}\label{sec 1}

In this section, we prove several lemmas necessary to prove the main results of the article.

We consider the set of all positive Lebesgue measurable functions $v(t)$ on $[1, \infty)$ for which the function $t^{-\varepsilon} v(t)$ decreases almost and the function $(\log \, 2t)^{\varepsilon} v (t) \uparrow$
  on $[1, \infty)$ for any number $\varepsilon> 0 $.  We will denote this set by $ SVL[1, \infty)$. It is clear that $SVL[1, \infty) \subset SV[1, \infty)$.

Next, we consider the sets
$Y^{m}(\bar{\gamma}, n)=\{\bar{s}\in\mathbb{Z}_{+}^{m}: \,\,
 \langle\bar{s},\bar{\gamma}\rangle\geq n \}$ and $\kappa^{m}(n, \bar{\gamma})=\{\bar{s}\in\mathbb{Z}_{+}^{m}: \,\, \langle\bar{s},\bar{\gamma}\rangle =n\}$, $n\in \mathbb{N}$. 
  
 \smallskip

\begin{lemma}\label{lem 1} {\it Let $\overline{\gamma}^{'} = (\gamma_{1}^{'}, \ldots, \gamma_{m}^{'})$,  $\overline{\gamma} = (\gamma_{1}, \ldots, \gamma_{m})$, $\overline{\theta}=(\theta_{1}, \ldots, \theta_{m})$ and the functions $v_{j}\in SLV[1, \infty)$, $0< \gamma_{j}^{'}\leq \gamma_{j}$, $1\leq \theta_{j}< \infty$, $j=1,\ldots , m$ and $\alpha \in (0, \infty)$, then 
$$
\left\|\left\{2^{-\alpha\langle\overline{s}, \overline{\gamma}\rangle}
\prod_{j=1}^{m}
V_{j}(2^{-s_{j}})
\right\}_{\bar{s}\in Y^{m}(n, \bar{\gamma}^{'})}\right\|_{
l_{\bar \theta}}\leq C2^{-n\alpha\delta}\prod_{j=1}^{m}
V_{j}(2^{-n})n^{\sum\limits_{j\in A\setminus\{j_{1}\}}\frac{1}{\theta_j}},
$$
for $n \in \mathbb{N}$, such that $n > n_{0}$, 
 where  $n_{0}$ is some positive number and  
 $\delta=\min\{\frac{\gamma_{j}}{\gamma_{j}^{'}} : j=1,\ldots, m\}$, $A=\{j : \frac{\gamma_{j}}{\gamma_{j}^{'}}=\delta, j=1,\ldots, m\}$, $j_{1}=\min\{j : j\in A\}$.
}
\end{lemma} 
{\bf Proof.} Let $m = 2$, then by the definition of the set
$Y^{2}(n, \bar{\gamma}^{'})$, we have that
$$
I_{n}:=\left\|\left\{2^{-\alpha\langle\overline{s}, \overline{\gamma}\rangle}
\prod_{j=1}^{2}
V_{j}(2^{-s_{j}})
\right\}_{\bar{s}\in Y^{2}(n, \bar{\gamma}^{'})}\right\|_{
l_{\bar \theta}} =
$$
$$
\Bigl\{\sum\limits_{s_{2}< \frac{n}{\gamma_{2}^{'}}}\Bigl[\sum\limits_{s_{1}\geq \frac{n-s_{2}\gamma_{2}^{'}}{\gamma_{1}^{'}}} 2^{-\alpha\langle\overline{s}, \overline{\gamma}\rangle \theta_{1}}\prod_{j=1}^{2}V_{j}^{\theta_{1}}(2^{-s_{j}}) \Bigr]^{\frac{\theta_{2}}{\theta_{1}}} + \sum\limits_{s_{2}\geq  \frac{n}{\gamma_{2}^{'}}}\Bigl[\sum\limits_{s_{1}=0}^{\infty} 2^{-\alpha\langle\overline{s}, \overline{\gamma}\rangle \theta_{1}}\prod_{j=1}^{2}V_{j}^{\theta_{1}}(2^{-s_{j}}) \Bigr]^{\frac{\theta_{2}}{\theta_{1}}} \Bigr\}^{\frac{1}{\theta_{2}}}.    \eqno (1) 
$$
Since the function $v_{1}\in SV[1, \infty)$, then the function $v_{1}(t)t^{-\varepsilon}$ almost decreases on $[1, \infty) $ for any numbers $\varepsilon> 0 $. Therefore, choosing the number $\varepsilon \in (0, \alpha)$, we will have
$$
\sum\limits_{s_{1}\geq \frac{n-s_{2}\gamma_{2}^{'}}{\gamma_{1}^{'}}} 2^{-s_{1}\gamma_{1}\alpha \theta_{1}}V_{1}^{\theta_{1}}(2^{-s_{1}}) \leq C\Bigl(v_{1}(2^{\frac{n-s_{2}\gamma_{2}^{'}}{\gamma_{1}^{'}}})(2^{\frac{n-s_{2}\gamma_{2}^{'}}{\gamma_{1}^{'}}})^{-\varepsilon} \Bigr)^{\theta_{1}}\sum\limits_{s_{1}\geq \frac{n-s_{2}\gamma_{2}^{'}}{\gamma_{1}^{'}}} 2^{-s_{1}\gamma_{1}\alpha \theta_{1}}2^{s_{1}\varepsilon\theta_{1}}
$$
$$
\leq CV_{1}^{\theta_{1}}(2^{-\frac{n-s_{2}\gamma_{2}^{'}}{\gamma_{1}^{'}}})2^{-\gamma_{1}\frac{n-s_{2}\gamma_{2}^{'}}{\gamma_{1}^{'}} \alpha \theta_{1}}. \eqno (2) 
$$
Since the function $v_{2}(t)t^{-\varepsilon}$ almost decreases on $[1, \infty)$, then arguing as in the proof of inequality (2), we can verify that
$$
\sum\limits_{s_{2}\geq \frac{n}{\gamma_{2}^{'}}} 2^{-s_{2}\gamma_{2}\alpha \theta_{2}}V_{1}^{\theta_{2}}(2^{-s_{2}}) \leq C V_{2}^{\theta_{2}}(2^{-\frac{n}{\gamma_{2}^{'}}})2^{-n\frac{\gamma_{2}}{\gamma_{2}^{'}}\alpha \theta_{2}}. \eqno  (3) 
$$
Now, it follows from inequalities (1), (2), and (3) that
$$
I_{n} \leq C \left\{2^{-n\frac{\gamma_{1}}{\gamma_{1}^{'}}\alpha \theta_{2}}\sum\limits_{s_{2}< \frac{n}{\gamma_{2}^{'}}} 2^{s_{2}(\frac{\gamma_{1}\gamma_{2}^{'}}{\gamma_{1}^{'}} - \gamma_{2})\alpha \theta_{2}}V_{2}^{\theta_{2}}(2^{-s_{2}}) V_{1}^{\theta_{2}}(2^{-\frac{n-s_{2}\gamma_{2}^{'}}{\gamma_{1}^{'}}})  \right.
$$
$$
\left. + 2^{-n\frac{\gamma_{2}}{\gamma_{2}^{'}}\alpha \theta_{2}}V_{2}^{\theta_{2}}(2^{-\frac{n}{\gamma_{2}^{'}}})\Bigl[\sum\limits_{s_{1}=0}^{\infty}2^{-s_{1}\gamma_{1}\alpha \theta_{1}}V_{1}^{\theta_{1}}(2^{-s_{1}}) \Bigr]^{\frac{\theta_{2}}{\theta_{1}}}
\right\}^{\frac{1}{\theta_{2}}}.  \eqno (4) 
$$
Further, it follows from inequality (2) that the series converges
$$
\sum\limits_{s_{1}=0}^{\infty}2^{-s_{1}\gamma_{1}\alpha \theta_{1}}V_{1}^{\theta_{1}}(2^{-s_{1}})
$$
converges. 
Therefore, inequality (4) implies that
$$
I_{n} \ll \left\{2^{-n\frac{\gamma_{1}}{\gamma_{1}^{'}}\alpha \theta_{2}}\sum\limits_{s_{2}< \frac{n}{\gamma_{2}^{'}}} 2^{-s_{2}\gamma_{2}^{'}(\frac{\gamma_{2}}{\gamma_{2}^{'}} -  \frac{\gamma_{1}}{\gamma_{1}^{'}})\alpha \theta_{2}}V_{2}^{\theta_{2}}(2^{-s_{2}}) V_{1}^{\theta_{2}}(2^{-\frac{n-s_{2}\gamma_{2}^{'}}{\gamma_{1}^{'}}}) \right.
$$
$$
\left.
  + 2^{-n\frac{\gamma_{2}}{\gamma_{2}^{'}}\alpha \theta_{2}}V_{2}^{\theta_{2}}(2^{-\frac{n}{\gamma_{2}^{'}}})
\right\}^{\frac{1}{\theta_{2}}}.  \eqno (5)
$$
Now, we estimate the sum in (5)

Let  
$\frac{\gamma_{2}}{\gamma_{2}^{'}} -  \frac{\gamma_{1}}{\gamma_{1}^{'}} > 0$. Choose the number $\lambda \in (0, (\frac{\gamma_{2}}{\gamma_{2}^{'}} -  \frac{\gamma_{1}}{\gamma_{1}^{'}})\alpha\gamma_{1}^{'})$. Then
$$
B_{n}=\sum\limits_{s_{2}< \frac{n}{\gamma_{2}^{'}}} 2^{-s_{2}\gamma_{2}^{'}(\frac{\gamma_{2}}{\gamma_{2}^{'}} -  \frac{\gamma_{1}}{\gamma_{1}^{'}})\alpha \theta_{2}}V_{2}^{\theta_{2}}(2^{-s_{2}}) V_{1}^{\theta_{2}}(2^{-\frac{n-s_{2}\gamma_{2}^{'}}{\gamma_{1}^{'}}}) =
$$ 
$$
=\sum\limits_{s_{2}< \frac{n}{\gamma_{2}^{'}}} 2^{-s_{2}\gamma_{2}^{'}(\frac{\gamma_{2}}{\gamma_{2}^{'}} -  \frac{\gamma_{1}}{\gamma_{1}^{'}})\alpha \theta_{2}}  2^{s_{2}\gamma_{2}^{'}\frac{\lambda}{\gamma_{1}^{'}}\theta_{2}}V_{2}^{\theta_{2}}(2^{-s_{2}}) 2^{-s_{2}\gamma_{2}^{'}\frac{\lambda}{\gamma_{1}^{'}}\theta_{2}}V_{1}^{\theta_{2}}(2^{-\frac{n-s_{2}\gamma_{2}^{'}}{\gamma_{1}^{'}}})
$$ 
 $$
= \sum\limits_{s_{2}< \frac{n}{\gamma_{2}^{'}}} 2^{-s_{2}\gamma_{2}^{'}((\frac{\gamma_{2}}{\gamma_{2}^{'}} -  \frac{\gamma_{1}}{\gamma_{1}^{'}})\alpha -\frac{\lambda}{\gamma_{1}^{'}}) \theta_{2}}V_{2}^{\theta_{2}}(2^{-s_{2}}) 2^{\frac{n-s_{2}\gamma_{2}^{'}}{\gamma_{1}^{'}}\lambda\theta_{2}}V_{1}^{\theta_{2}}(2^{-\frac{n-s_{2}\gamma_{2}^{'}}{\gamma_{1}^{'}}})2^{-n\frac{\lambda}{\gamma_{1}^{'}}\theta_{2}}. \eqno (6)
 $$
 Since $0< n-s_{2}\gamma_{2}^{'} \leq n$ and the function $v_{1}(t)t^{\varepsilon}$ is almost increasing by $[1, \infty)$, then
$$
2^{\frac{n-s_{2}\gamma_{2}^{'}}{\gamma_{1}^{'}}\lambda\theta_{2}}V_{1}^{\theta_{2}}(2^{-\frac{n-s_{2}\gamma_{2}^{'}}{\gamma_{1}^{'}}})\leq C2^{\frac{n}{\gamma_{1}^{'}}\lambda\theta_{2}} V_{1}^{\theta_{2}}(2^{-\frac{n}{\gamma_{1}^{'}}}).
$$
Hence, from equalities (6), we obtain
$$
B_{n}\leq CV_{1}^{\theta_{2}}(2^{-\frac{n}{\gamma_{1}^{'}}})\sum\limits_{s_{2}< \frac{n}{\gamma_{2}^{'}}} 2^{-s_{2}\gamma_{2}^{'}((\frac{\gamma_{2}}{\gamma_{2}^{'}} -  \frac{\gamma_{1}}{\gamma_{1}^{'}})\alpha -\frac{\lambda}{\gamma_{1}^{'}}) \theta_{2}}V_{2}^{\theta_{2}}(2^{-s_{2}}). \eqno (7)
$$
Since $(\frac{\gamma_{2}}{\gamma_{2}^{'}} -  \frac{\gamma_{1}}{\gamma_{1}^{'}})\alpha -\frac{\lambda}{\gamma_{1}^{'}} > 0$ and the function $v_{2}\in SV[1, \infty)$, then
$$
\sum\limits_{s_{2}=0}^{\infty} 2^{-s_{2}\gamma_{2}^{'}((\frac{\gamma_{2}}{\gamma_{2}^{'}} -  \frac{\gamma_{1}}{\gamma_{1}^{'}})\alpha -\frac{\lambda}{\gamma_{1}^{'}})\theta_{2}}V_{2}^{\theta_{2}}(2^{-s_{2}})< \infty.
$$
Therefore, from inequality (7), it follows that  
$$
B_{n}\leq CV_{1}^{\theta_{2}}(2^{-\frac{n}{\gamma_{1}^{'}}})  \eqno (8) 
$$
in the case 
 $\frac{\gamma_{2}}{\gamma_{2}^{'}} -  \frac{\gamma_{1}}{\gamma_{1}^{'}} > 0$.

Let 
 $\frac{\gamma_{2}}{\gamma_{2}^{'}} -  \frac{\gamma_{1}}{\gamma_{1}^{'}} < 0$, then 
  $\frac{\gamma_{1}}{\gamma_{1}^{'}} -  \frac{\gamma_{2}}{\gamma_{2}^{'}} > 0$. 
We choose the number $\eta\in (0, (\frac{\gamma_{1}}{\gamma_{1}^{'}} -  \frac{\gamma_{2}}{\gamma_{2}^{'}})\alpha)$,
then  
  $$
B_{n}=\sum\limits_{s_{2}< \frac{n}{\gamma_{2}^{'}}} 2^{-s_{2}\gamma_{2}^{'}(\frac{\gamma_{2}}{\gamma_{2}^{'}} -  \frac{\gamma_{1}}{\gamma_{1}^{'}})\alpha \theta_{2}}V_{2}^{\theta_{2}}(2^{-s_{2}}) V_{1}^{\theta_{2}}(2^{-\frac{n-s_{2}\gamma_{2}^{'}}{\gamma_{1}^{'}}}) =
$$ 
$$
=\sum\limits_{s_{2}< \frac{n}{\gamma_{2}^{'}}} 2^{-s_{2}\gamma_{2}^{'}(\frac{\gamma_{2}}{\gamma_{2}^{'}} -  \frac{\gamma_{1}}{\gamma_{1}^{'}})\alpha \theta_{2}}  2^{s_{2}\gamma_{2}^{'}\eta\theta_{2}}V_{2}^{\theta_{2}}(2^{-s_{2}}) 2^{-s_{2}\gamma_{2}^{'}\eta\theta_{2}}V_{1}^{\theta_{2}}(2^{-\frac{n-s_{2}\gamma_{2}^{'}}{\gamma_{1}^{'}}})
$$ 
 $$
= \sum\limits_{s_{2}< \frac{n}{\gamma_{2}^{'}}} 2^{-s_{2}\gamma_{2}^{'}((\frac{\gamma_{2}}{\gamma_{2}^{'}} -  \frac{\gamma_{1}}{\gamma_{1}^{'}})\alpha +\eta) \theta_{2}}V_{2}^{\theta_{2}}(2^{-s_{2}}) 2^{s_{2}\gamma_{2}^{'}\eta\theta_{2}}V_{1}^{\theta_{2}}(2^{-\frac{n-s_{2}\gamma_{2}^{'}}{\gamma_{1}^{'}}}). \eqno (9)  
 $$
 Since the function $v_{2}(t)t^{\varepsilon}$ almost increases by $[1, \infty)$, then
 $$
 2^{s_{2}\gamma_{2}^{'}\eta\theta_{2}}V_{2}^{\theta_{2}}(2^{-s_{2}})\leq C 2^{n\eta\theta_{2}}V_{2}^{\theta_{2}}(2^{-\frac{n}{\gamma_{2}^{'}}})
 $$
for  $0\leq s_{2}< \frac{n}{\gamma_{2}^{'}}$.
  Therefore, from equalities (9), we obtain
  $$
B_{n}\leq C2^{n\eta\theta_{2}}V_{2}^{\theta_{2}}(2^{-\frac{n}{\gamma_{2}^{'}}})\sum\limits_{s_{2}< \frac{n}{\gamma_{2}^{'}}} 2^{-s_{2}\gamma_{2}^{'}((\frac{\gamma_{2}}{\gamma_{2}^{'}} -  \frac{\gamma_{1}}{\gamma_{1}^{'}})\alpha +\eta) \theta_{2}}V_{1}^{\theta_{2}}(2^{-\frac{n-s_{2}\gamma_{2}^{'}}{\gamma_{1}^{'}}})
$$
$$
\leq C2^{n\eta\theta_{2}}V_{2}^{\theta_{2}}(2^{-\frac{n}{\gamma_{2}^{'}}})\sum\limits_{s_{2}< \frac{n}{\gamma_{2}^{'}}} 2^{(n-s_{2}\gamma_{2}^{'})((\frac{\gamma_{2}}{\gamma_{2}^{'}} -  \frac{\gamma_{1}}{\gamma_{1}^{'}})\alpha +\eta) \theta_{2}}V_{1}^{\theta_{2}}(2^{-\frac{n-s_{2}\gamma_{2}^{'}}{\gamma_{1}^{'}}})2^{-n((\frac{\gamma_{2}}{\gamma_{2}^{'}} -  \frac{\gamma_{1}}{\gamma_{1}^{'}})\alpha +\eta) \theta_{2}}. \eqno (10) 
$$
Since $(\frac{\gamma_{2}}{\gamma_{2}^{'}} -  \frac{\gamma_{1}}{\gamma_{1}^{'}})\alpha -\eta < 0$ and the function $v_{1} \in SV[1, \infty)$, then
 $$
\sum\limits_{s_{2}< \frac{n}{\gamma_{2}^{'}}} 2^{(n-s_{2}\gamma_{2}^{'})((\frac{\gamma_{2}}{\gamma_{2}^{'}} -  \frac{\gamma_{1}}{\gamma_{1}^{'}})\alpha +\eta) \theta_{2}}V_{1}^{\theta_{2}}(2^{-\frac{n-s_{2}\gamma_{2}^{'}}{\gamma_{1}^{'}}})
=\sum\limits_{0<l\leq \frac{n}{\gamma_{1}^{'}}} 2^{-l((\frac{\gamma_{1}}{\gamma_{1}^{'}} -  \frac{\gamma_{2}}{\gamma_{2}^{'}})\alpha -\eta) \theta_{2}}V_{1}^{\theta_{2}}(2^{-l})
$$
$$
\leq \sum\limits_{s_{2}=0}^{\infty} 2^{-l((\frac{\gamma_{1}}{\gamma_{1}^{'}} -  \frac{\gamma_{2}}{\gamma_{2}^{'}})\alpha -\eta) \theta_{2}}V_{1}^{\theta_{2}}(2^{-l})< \infty.
$$
Therefore,  from inequality (10), it follows that
$$
B_{n}\leq C2^{-n((\frac{\gamma_{2}}{\gamma_{2}^{'}} -  \frac{\gamma_{1}}{\gamma_{1}^{'}})\alpha ) \theta_{2}}V_{2}^{\theta_{2}}(2^{-\frac{n}{\gamma_{2}^{'}}})  \eqno (11)  
$$
in the case 
 $\frac{\gamma_{2}}{\gamma_{2}^{'}} -  \frac{\gamma_{1}}{\gamma_{1}^{'}} < 0$.

Let 
 $\frac{\gamma_{2}}{\gamma_{2}^{'}} -  \frac{\gamma_{1}}{\gamma_{1}^{'}} = 0$, then 
$$
 \sum\limits_{0\leq s_{2}< \frac{n}{\gamma_{2}^{'}}} V_{2}^{\theta_{2}}(2^{-s_{2}}) V_{1}^{\theta_{2}}(2^{-\frac{n-s_{2}\gamma_{2}^{'}}{\gamma_{1}^{'}}}) + V_{2}^{\theta_{2}}(2^{-\frac{n}{\gamma_{2}^{'}}}) = 
 $$
$$
 \sum\limits_{0\leq s_{2}< \frac{n}{\gamma_{2}^{'}}} V_{2}^{\theta_{2}}(2^{-s_{2}}) V_{1}^{\theta_{2}}(2^{-\frac{n-s_{2}\gamma_{2}^{'}}{\gamma_{1}^{'}}}) + V_{1}^{\theta_{2}}(1)V_{2}^{\theta_{2}}(2^{-\frac{n}{\gamma_{1}^{'}}})\frac{1}{V_{1}^{\theta_{2}}(1)}   
$$
$$ 
 \leq \max\{1, \frac{1}{V_{1}^{\theta_{2}}(1)}\} \sum\limits_{0\leq s_{2}< \frac{n}{\gamma_{2}^{'}}} V_{2}^{\theta_{2}}(2^{-s_{2}}) V_{1}^{\theta_{2}}(2^{-\frac{n-s_{2}\gamma_{2}^{'}}{\gamma_{1}^{'}}}). \eqno (12) 
$$
Since the functions $v_{j}\in SVL[1, \infty )$, $j=1, 2$, then the functions $v_{j}(t)(\log 2t)^{\varepsilon}$, $j = 1, 2$ increase on $ [1, \infty) $ for $\varepsilon> 0$.
Therefore
$$
V_{2}^{\theta_{2}}(2^{-s_{2}})=v_{2}^{\theta_{2}}(2^{s_{2}})(\log 2^{s_{2}+1})^{\varepsilon}(\log 2^{s_{2}+1})^{-\varepsilon}\leq v_{2}^{\theta_{2}}(2^{\frac{n}{\gamma_{2}^{'}}})\Bigl(\log 2^{\frac{n}{\gamma_{2}^{'}}+1}\Bigr)^{\varepsilon}(\log 2^{s_{2}+1})^{-\varepsilon} 
$$
$$
= V_{2}^{\theta_{2}}(2^{-\frac{n}{\gamma_{2}^{'}}})\Bigl(\frac{n}{\gamma_{2}^{'}}+1\Bigr)^{\varepsilon\theta_{2}}(s_{2}+1)^{-\varepsilon\theta_{2}} 
$$
and
$$
V_{1}^{\theta_{2}}(2^{-\frac{n-s_{2}\gamma_{2}^{'}}{\gamma_{1}^{'}}}) = v_{1}^{\theta_{2}}(2^{\frac{n-s_{2}\gamma_{2}^{'}}{\gamma_{1}^{'}}}) = \Bigl(\log 2^{1+\frac{n-s_{2}\gamma_{2}^{'}}{\gamma_{1}^{'}}}\Bigr)^{\varepsilon\theta_{2}}v_{1}^{\theta_{2}}(2^{\frac{n-s_{2}\gamma_{2}^{'}}{\gamma_{1}^{'}}})\Bigl(\log 2^{1+\frac{n-s_{2}\gamma_{2}^{'}}{\gamma_{1}^{'}}}\Bigr)^{-\varepsilon\theta_{2}} 
$$
$$
\leq C
\Bigl(1+\frac{n}{\gamma_{1}^{'}}\Bigr)^{\varepsilon\theta_{2}}v_{1}^{\theta_{2}}(2^{\frac{n}{\gamma_{1}^{'}}})\Bigl(1+\frac{n-s_{2}\gamma_{2}^{'}}{\gamma_{1}^{'}}\Bigr)^{-\varepsilon\theta_{2}} 
$$
$$
=  C
\Bigl(1+\frac{n}{\gamma_{1}^{'}}\Bigr)^{\varepsilon\theta_{2}}V_{1}^{\theta_{2}}(2^{-\frac{n}{\gamma_{1}^{'}}})\Bigl(1+\frac{n-s_{2}\gamma_{2}^{'}}{\gamma_{1}^{'}}\Bigr)^{-\varepsilon\theta_{2}}. 
$$
Hence 
$$
\sum\limits_{0\leq s_{2}< \frac{n}{\gamma_{2}^{'}}} V_{2}^{\theta_{2}}(2^{-s_{2}}) V_{1}^{\theta_{2}}(2^{-\frac{n-s_{2}\gamma_{2}^{'}}{\gamma_{1}^{'}}})
$$
$$
\leq C\Bigl(1+\frac{n}{\gamma_{1}^{'}}\Bigr)^{\varepsilon\theta_{2}}V_{1}^{\theta_{2}}(2^{-\frac{n}{\gamma_{1}^{'}}})V_{2}^{\theta_{2}}(2^{-\frac{n}{\gamma_{2}^{'}}})\Bigl(\frac{n}{\gamma_{2}^{'}}+1\Bigr)^{\varepsilon\theta_{2}}
\sum\limits_{0\leq s_{2}< \frac{n}{\gamma_{2}^{'}}}\Bigl(1+\frac{n-s_{2}\gamma_{2}^{'}}{\gamma_{1}^{'}}\Bigr)^{-\varepsilon\theta_{2}}(s_{2}+1)^{-\varepsilon\theta_{2}}.  \eqno (13)
$$
Now, choosing the number $\varepsilon \in (0, \frac{1}{\theta_{2}})$ and using Lemma 1 [19] from inequality (13), we obtain
$$
\sum\limits_{0\leq s_{2}< \frac{n}{\gamma_{2}^{'}}} V_{2}^{\theta_{2}}(2^{-s_{2}}) V_{1}^{\theta_{2}}(2^{-\frac{n-s_{2}\gamma_{2}^{'}}{\gamma_{1}^{'}}})
\leq CV_{1}^{\theta_{2}}(2^{-\frac{n}{\gamma_{1}^{'}}})V_{2}^{\theta_{2}}(2^{-\frac{n}{\gamma_{2}^{'}}})(n+1) \eqno ((14)
$$
in the case 
$\frac{\gamma_{2}}{\gamma_{2}^{'}} -  \frac{\gamma_{1}}{\gamma_{1}^{'}} = 0$.

Further, it follows from inequalities (5), (8), (11), and (14) that
$$
I_{n} \leq C \left\{2^{-n\frac{\gamma_{1}}{\gamma_{1}^{'}}\alpha}V_{1}^{\theta_{2}}(2^{-\frac{n}{\gamma_{1}^{'}}}) + 2^{-n\frac{\gamma_{2}}{\gamma_{2}^{'}}\alpha \theta_{2}}V_{2}^{\theta_{2}}(2^{-\frac{n}{\gamma_{2}^{'}}})\right\}^{\frac{1}{\theta_{2}}} \eqno (15)  
$$
in the case 
 $\frac{\gamma_{2}}{\gamma_{2}^{'}} -  \frac{\gamma_{1}}{\gamma_{1}^{'}} > 0$ and
$$
I_{n} \leq C \left\{2^{-n\frac{\gamma_{2}}{\gamma_{2}^{'}}\alpha}V_{2}^{\theta_{2}}(2^{-\frac{n}{2}}) + 2^{-n\frac{\gamma_{2}}{\gamma_{2}^{'}}\alpha \theta_{2}}V_{2}^{\theta_{2}}(2^{-n})\right\}^{\frac{1}{\theta_{2}}} \eqno (16)    
$$
in the case 
 $\frac{\gamma_{2}}{\gamma_{2}^{'}} -  \frac{\gamma_{1}}{\gamma_{1}^{'}} < 0$ and

$$
I_{n} \leq  C2^{-n\frac{\gamma_{1}}{\gamma_{1}^{'}}\alpha}V_{1}(2^{-\frac{n}{2}})V_{2}(2^{-\frac{n}{2}})(n+1)^{\frac{1}{\theta_{2}}} \eqno (17)    
$$
in the case 
 $\frac{\gamma_{2}}{\gamma_{2}^{'}} -  \frac{\gamma_{1}}{\gamma_{1}^{'}} = 0$.

Since the functions $v_{j}\in SVL[1, \infty )\subset  SV[1, \infty )$, $ j = 1,2 $, then the function $\frac{v_{2}}{v_{1}} \in SV[1, \infty)$ (see \cite[p. 16]{2}).
It is known that if the function $w\in SV[1, \infty)$, then $x^{-\varepsilon}w(x)\rightarrow 0$ for $ x \rightarrow \infty$, for any number $\varepsilon > 0$ (see \cite[p. 52]{1}, \cite[p. 109]{3}).
Therefore   
$$
\frac{2^{-n\gamma_{2}^{'}(\frac{\gamma_{2}}{\gamma_{2}^{'}} -  \frac{\gamma_{1}}{\gamma_{1}^{'}})\alpha}V_{2}(2^{-n})}{V_{1}(2^{-n})} = \frac{2^{-n\gamma_{2}^{'}(\frac{\gamma_{2}}{\gamma_{2}^{'}} -  \frac{\gamma_{1}}{\gamma_{1}^{'}})\alpha}v_{2}(2^{n})}{v_{1}(2^{n})} \rightarrow 0
$$
for $n \rightarrow \infty$, in the case 
 $\frac{\gamma_{2}}{\gamma_{2}^{'}} -  \frac{\gamma_{1}}{\gamma_{1}^{'}} > 0$. 
Hence 
$$
2^{-n\frac{\gamma_{2}}{\gamma_{2}^{'}}\alpha}V_{2}(2^{-n})=2^{-n\frac{\gamma_{1}}{\gamma_{1}^{'}}\alpha}V_{1}(2^{-n})V_{2}(2^{-n}) \frac{2^{-n\gamma_{2}^{'}(\frac{\gamma_{2}}{\gamma_{2}^{'}} -  \frac{\gamma_{1}}{\gamma_{1}^{'}})\alpha}}{V_{1}(2^{-n})}  \leq 2^{-n\frac{\gamma_{1}}{\gamma_{1}^{'}}\alpha}V_{1}(2^{-n})   \eqno (18)    
$$
for natural numbers $n> n_{0} $, where $n_{0}$ is some positive number, in the case of $\frac{\gamma_{2}}{\gamma_{2}^{'}} -  \frac{\gamma_{1}}{\gamma_{1}^{'}} > 0$.
  Now, it follows from inequalities (15) and (18) that
 $$
I_{n} \leq C2^{-n\frac{\gamma_{1}}{\gamma_{1}^{'}}\alpha}V_{1}(2^{-\frac{n}{\gamma_{1}^{'}}})   \eqno   (19) 
$$
for natural numbers $ n> n_{0}$, in the case of $\frac{\gamma_{2}}{\gamma_{2}^{'}} -  \frac{\gamma_{1}}{\gamma_{1}^{'}} > 0$.

Since the functions $v_{j}\in SVL[1, \infty )$, $j = 1, 2$, then $V_{j}(2^{-\frac{n}{\gamma_{j}^{'}}})\ll V_{j}(2^{-n})$. Therefore,
estimates (16), (17), and (19) can be written in the following form
$$
I_{n} \leq C2^{-n\alpha\delta}\prod_{j\in A}^{}V_{j}(2^{-n})(n+1)^{\sum\limits_{j\in A}\frac{1}{\theta_{j}}},   \eqno (20)  
$$
where $\delta=\min\{\frac{\gamma_{j}}{\gamma_{j}^{'}}: \,\, j=1, 2\}$ and $A=\{j : \delta=\frac{\gamma_{j}}{\gamma_{j}^{'}}, j=1, 2\}$.
 This proves the lemma for $ m = 2 $.

Now, we will assume that the lemma is true for $m-1 \geq 2$, that is
$$
\left\|\left\{2^{-\alpha\langle\overline{s}_{m-1}, \overline{\gamma}_{m-1}\rangle}
\prod_{j=1}^{m-1}
V_{j}(2^{-s_{j}})
\right\}_{\bar{s}\in Y^{m-1}(n, \bar{\gamma}_{m-1}^{'})}\right\|_{
l_{\bar \theta}(\mathbb{Z}^{m-1}_{+})}
$$
$$
\leq C2^{-n\alpha\delta_{m-1}}\prod_{j\in A_{m-1}}
V_{j}(2^{-n})n^{\sum\limits_{j\in A_{m-1}\setminus\{j_{0}\}}\frac{1}{\theta_j}},
\eqno (21)  
$$
where 
 $\delta_{m-1}=\min\{\frac{\gamma_{j}}{\gamma_{j}^{'}} : j=1,\ldots, m-1\}$, $A_{m-1}=\{j : \frac{\gamma_{j}}{\gamma_{j}^{'}}=\delta, j=1,\ldots, m-1\}$, $j_{0}=\min\{j : j\in A_{m-1}\}$, $\overline{a}_{m-1}=(a_{1},...,a_{m-1})$.

Let us prove the lemma for $ m $. By the definition of the set $Y^{m}(n, \bar{\gamma}^{'})$ and by the Minkowski inequality, we obtain
$$
\left\|\left\{2^{-\alpha\langle\overline{s}, \overline{\gamma}\rangle}
\prod_{j=1}^{m}
V_{j}(2^{-s_{j}})
\right\}_{\bar{s}\in Y^{m}(n, \bar{\gamma}^{'})}\right\|_{
l_{\bar \theta}}\ll 
$$
$$
 \left\{\left[\sum\limits_{0\leq s_{m}< \frac{n}{\gamma_{m}^{'}}}\Bigl(2^{-s_{m}\gamma_{m}\alpha}V_{m}(2^{-s_{m}}) 
J_{n-s_{m}\gamma_{m}}\Bigr)^{\theta_{m}}\right]^{\frac{1}{\theta_{m}}} \right. 
$$
$$
\left. +  \left[\sum\limits_{ s_{m}\geq  \frac{n}{\gamma_{m}^{'}}}\Bigl(2^{-s_{m}\gamma_{m}\alpha}V_{m}(2^{-s_{m}})\left\|\left\{2^{-\alpha\langle\overline{s}_{m-1}, \overline{\gamma}_{m-1}\rangle}
\prod_{j=1}^{m-1}
V_{j}(2^{-s_{j}})
\right\}_{\bar{s}_{m-1}\in \mathbb{Z}_{+}^{m-1}}\right\|_{
l_{\bar \theta}(\mathbb{Z}_{+}^{m-1})} \Bigr)^{\theta_{m}} \right]^{\frac{1}{\theta_{m}}}\right\} 
$$
$$
= C\left\{\sigma_{1}(n) + \sigma_{2}(n)\right\},  \eqno (22)    
$$
where 
$$
J_{n-s_{m}\gamma_{m}}=\left\|\left\{2^{-\alpha\langle\overline{s}_{m-1}, \overline{\gamma}_{m-1}\rangle}
\prod_{j=1}^{m-1}
V_{j}(2^{-s_{j}})
\right\}_{\bar{s}\in Y^{m-1}(n-s_{m}\gamma_{m}, \bar{\gamma}_{m-1}^{'})}\right\|_{
l_{\bar \theta}(\mathbb{Z}_{+}^{m-1})}.
$$
To estimate $\sigma_{1}(n)$, we use inequality (21), then we get
$$
\sigma_{1}(n) \leq C\left[\sum\limits_{0\leq s_{m}< \frac{n}{\gamma_{m}^{'}}}\Bigl(2^{-s_{m}\gamma_{m}\alpha}V_{m}(2^{-s_{m}}) \right.
$$
$$
\left.
\times 2^{-(n-s_{m}\gamma_{m}^{'})\alpha\delta_{m-1}}\prod_{j\in A_{m-1}}
V_{j}(2^{-(n-s_{m}\gamma_{m}^{'})})(n-s_{m}\gamma_{m}^{'})^{\sum\limits_{j\in A_{m-1}\setminus\{j_{0}\}}\frac{1}{\theta_j}} \Bigr)^{\theta_{m}} \right]^{\frac{1}{\theta_{m}}}
$$
$$
\leq C2^{-n\alpha\delta_{m-1}}\left[\sum\limits_{0\leq s_{m}< \frac{n}{\gamma_{m}^{'}}}2^{-s_{m}\gamma_{m}^{'}(\frac{\gamma_{m}}{\gamma_{m}^{'}}-\delta_{m-1})\alpha\theta_{m}}V_{m}^{\theta_{m}}(2^{-s_{m}}) \right.
$$
$$
\left.
\prod_{j\in A_{m-1}}
V_{j}^{\theta_{m}}(2^{-(n-s_{m}\gamma_{m}^{'})})(n-s_{m}\gamma_{m}^{'})^{\theta_{m}\sum\limits_{j\in A_{m-1}\setminus\{j_{0}\}}\frac{1}{\theta_j}}\right]^{\frac{1}{\theta_{m}}}.  \eqno (23)  
$$

If 
 $\frac{\gamma_{m}}{\gamma_{m}^{'}}-\delta_{m-1} >0$, then  
$\delta_{m}=\min\{\frac{\gamma_{j}}{\gamma_{j}^{'}} : j=1,\ldots, m\}=\delta_{m-1}$, $A_{m}=\{j : \frac{\gamma_{j}}{\gamma_{j}^{'}}=\delta, j=1,\ldots, m\}=A_{m-1}$, $j_{1}=\min\{j : j\in A_{m}\}=j_{0}$, $j^{'}=\max\{j\in A_{m-1}\}$, $\overline{a}_{m-1}=(a_{1},...,a_{m-1})$.
Therefore 
$$
\sum\limits_{j\in A_{m}\setminus\{j_{1}\}}\frac{1}{\theta_j}=\sum\limits_{j\in A_{m-1}\setminus\{j_{0}\}}\frac{1}{\theta_j}.
$$ 
Now, taking into account that the functions $v_{j}\in SVL[1, \infty )$, $ j = 1, \ldots, m $ and arguing as in the proof of inequality (8) from (23), we obtain
$$
\sigma_{1}(n) \leq C2^{-n\alpha\delta_{m-1}}\prod_{j\in A_{m-1}}
V_{j}(2^{-n})(n+1)^{\sum\limits_{j\in A_{m-1}\setminus\{j_{0}\}}\frac{1}{\theta_j}}
$$
$$
=C2^{-n\alpha\delta_{m}}\prod_{j\in A_{m}}
V_{j}(2^{-n})(n+1)^{\sum\limits_{j\in A_{m}\setminus\{j_{0}\}}\frac{1}{\theta_j}}. \eqno (24)  
$$
Let 
 $\frac{\gamma_{m}}{\gamma_{m}^{'}}-\delta_{m-1} <0$, then 
 $\frac{\gamma_{m}}{\gamma_{m}^{'}}< \frac{\gamma_{j}}{\gamma_{j}^{'}}$, $j=1,...,m-1$. 
Since
$\delta_{m-1}-\frac{\gamma_{m}}{\gamma_{m}^{'}} > 0$, then taking into account that the functions $v_{j} \in SVL[1, \infty)$, $ j = 1, \ldots, m$ and arguing as in the proof of inequality (11) from (23), we obtain
$$
\sigma_{1}(n) \leq C2^{-n\alpha\delta_{m-1}}2^{-n(\frac{\gamma_{m}}{\gamma_{m}^{'}}-\delta_{m-1})\alpha}V_{m}(2^{-n})=C2^{-n\alpha\delta_{m}}V_{m}(2^{-n}). \eqno (25)   
$$
Let $\frac{\gamma_{m}}{\gamma_{m}^{'}}-\delta_{m-1} =0$, then 
$A_{m-1}\subset A_{m}$ and $j_{1}=\min\{j\in A_{m}\}=\min\{j\in A_{m-1}\}=j_{0}$.
Therefore, from inequality (23), we obtain 
$$
\sigma_{1}(n) \leq C2^{-n\alpha\delta_{m-1}}
$$
$$
\times \left[\sum\limits_{0\leq s_{m}< \frac{n}{\gamma_{m}^{'}}}V_{m}^{\theta_{m}}(2^{-s_{m}})
\prod_{j\in A_{m-1}}
V_{j}^{\theta_{m}}(2^{-(n-s_{m}\gamma_{m}^{'})})(n-s_{m}\gamma_{m}^{'})^{\theta_{m}\sum\limits_{j\in A_{m-1}\setminus\{j_{0}\}}\frac{1}{\theta_j}}\right]^{\frac{1}{\theta_{m}}}.  \eqno (26)  
$$
Further, taking into account that the functions $v_{j}\in SVL[1, \infty )$, $ j = 1, \ldots, m $ and arguing as in the proof of inequality (14) from (26), we obtain
$$
\sigma_{1}(n) \leq C2^{-n\alpha\delta_{m}}\prod_{j\in A_{m}}
V_{j}(2^{-n/2})(n-s_{m}\gamma_{m}^{'})^{\sum\limits_{j\in A_{m}\setminus\{j_{1}\}}\frac{1}{\theta_j}}  \eqno (27)  
$$
in the case 
 $\frac{\gamma_{m}}{\gamma_{m}^{'}}-\delta_{m-1} =0$.

Since (see (2) and (3)) 
$$
\left\|\left\{2^{-\alpha\langle\overline{s}_{m-1}, \overline{\gamma}_{m-1}\rangle}
\prod_{j=1}^{m-1}V_{j}(2^{-s_{j}})
\right\}_{\bar{s}_{m-1}\in \mathbb{Z}_{+}^{m-1}}\right\|_{
l_{\bar \theta}(\mathbb{Z}_{+}^{m-1})} < \infty
$$
and
$$
\sum\limits_{ s_{m}\geq  \frac{n}{\gamma_{m}^{'}}}\Bigl(2^{-s_{m}\gamma_{m}\alpha}V_{m}(2^{-s_{m}})\Bigr)^{\theta_{m}} \leq C\Bigl(2^{-n\frac{\gamma_{m}}{\gamma_{m}^{'}}\alpha}V_{m}(2^{-n})\Bigr)^{\theta_{m}},
$$
then 
$$
\sigma_{2}(n) \leq C2^{-n\frac{\gamma_{m}}{\gamma_{m}^{'}}\alpha}V_{m}(2^{-n})\leq C2^{-n\alpha\delta_{m}}\prod_{j\in A_{m}}
V_{j}(2^{-n})n^{\sum\limits_{j\in A_{m}\setminus\{j_{1}\}}\frac{1}{\theta_j}},
\eqno (28)   
$$
in the case  
$\frac{\gamma_{m}}{\gamma_{m}^{'}}-\delta_{m-1} \neq 0$.

Now, from inequalities (22), (24), (25), (27), (28), it follows that
$$
\left\|\left\{2^{-\alpha\langle\overline{s}, \overline{\gamma}\rangle}
\prod_{j=1}^{m}
V_{j}(2^{-s_{j}})
\right\}_{\bar{s}\in Y^{m}(n, \bar{\gamma}^{'})}\right\|_{
l_{\bar \theta}}\leq C2^{-n\alpha\delta_{m}}\prod_{j\in A_{m}}
V_{j}(2^{-n})n^{\sum\limits_{j\in A_{m}\setminus\{j_{1}\}}\frac{1}{\theta_j}}.
$$
If $A_{m}\setminus\{j_{1}\} = \emptyset$, then it is considered that $\sum\limits_{j\in A_{m}\setminus\{j_{1}\}}\frac{1}{\theta_j}=0 $.
Lemma 1 is proved.

\begin{lemma}\label{lem 2} 
{\it Let 
   $\overline{\gamma} = (\gamma_{1}, \ldots, \gamma_{m})$, $\overline{\varepsilon}=(\varepsilon_{1}, \ldots, \varepsilon_{m})$ and  $0<  \gamma_{j}$, $0< \varepsilon_{j}\leq \infty$, $j=1,\ldots, m$, $m\geq 2$, $\alpha \in (0, \infty)$ and functions    $v_{j}\in SVL[1, \infty )$, $j=1,\ldots, m$, then  
$$
J_{n}:=\left\|\left\{2^{-\alpha\langle\overline{s}, \overline{\gamma}\rangle}
\prod_{j=1}^{m}
V_{j}(2^{-s_{j}})
\right\}_{\bar{s}\in \kappa^{m}(n, \bar{\gamma})}\right\|_{
l_{\bar\varepsilon}}\gg 2^{-n\alpha}\prod_{j=1}^{m}
V_{j}(2^{-n})n^{\sum\limits_{j=2}^{m}\frac{1}{\varepsilon_j}}.
$$
}
\end{lemma}   

{\bf Proof.} Let $m =2$, then, by the definition of the set
$\kappa^{2}(n, \bar{\gamma}^{'})$ and considering that the functions $v_{j}\in SVL[1, \infty )$, $j = 1, 2$,  we have
$$
J_{n}:=\left\|\left\{2^{-\alpha\langle\overline{s}, \overline{\gamma}\rangle}
\prod_{j=1}^{2}V_{j}(2^{-s_{j}})
\right\}_{\bar{s}\in \kappa^{2}(n, \bar{\gamma})}\right\|_{
l_{\bar \varepsilon}} \geq
$$
$$
\Bigl\{\sum\limits_{s_{2}< \frac{n}{\gamma_{2}}}\Bigl[\sum\limits_{s_{1}= \frac{n-s_{2}\gamma_{2}}{\gamma_{1}}} 2^{-\alpha\langle\overline{s}, \overline{\gamma}\rangle \varepsilon_{1}}\prod_{j=1}^{2}V_{j}^{\varepsilon_{1}}(2^{-s_{j}}) \Bigr]^{\frac{\varepsilon_{2}}{\varepsilon_{1}}}  \Bigr\}^{\frac{1}{\varepsilon_{2}}}
$$
$$
= 2^{-n\alpha}\Bigl\{\sum\limits_{s_{2}< \frac{n}{\gamma_{2}}}V_{2}^{\varepsilon_{2}}(2^{-s_{2}})V_{1}^{\varepsilon_{2}}(2^{-(n - s_{2}\gamma_{2})}) \Bigr\}^{\frac{1}{\varepsilon_{2}}}\gg 2^{-n\alpha}\prod_{j=1}^{2}V_{j}(2^{-n})(n+1)^{\frac{1}{\varepsilon_{2}}}.    \eqno (29)    
$$
Now, we will assume that the lemma is true for $m-1 \geq 2$, that is
$$
\left\|\left\{2^{-\alpha\langle\overline{s}_{m-1}, \overline{\gamma}_{m-1}\rangle}
\prod_{j=1}^{m-1}
V_{j}(2^{-s_{j}})
\right\}_{\bar{s}_{m-1}\in \kappa^{m-1}(n, \bar{\gamma})}\right\|_{
l_{\bar\varepsilon}(\mathbb{Z}^{m-1}_{+})}
$$
$$
\gg 2^{-n\alpha}\prod_{j=1}^{m-1}V_{j}(2^{-n})(n+1)^{\sum\limits_{j=2}^{m-1}\frac{1}{\varepsilon_{j}}},
\eqno (30)   
$$
Let us prove the assertion of the lemma for $m$.
By the definition of the set
$\kappa^{m}(n, \bar{\gamma}^{'})$, according to assumption (30) and taking into account that the functions $v_{j}\in SVL[1, \infty )$, $j = 1, \ldots, m$, we have  
$$
J_{n} = 
$$
$$
\Biggl\{\sum\limits_{s_{m}< \frac{n}{\gamma_{m}}}2^{-s_{m}\gamma_{m}\alpha\varepsilon_{m}}V_{m}^{\varepsilon_{m}}(2^{-s_{m}})\left\|\left\{2^{-\alpha\langle\overline{s}_{m-1}, \overline{\gamma}_{m-1}\rangle}
\prod_{j=1}^{m-1}
V_{j}(2^{-s_{j}})
\right\}_{\bar{s}_{m-1}\in \kappa^{m-1}(n, \bar{\gamma})}\right\|_{
l_{\bar \varepsilon}(\mathbb{Z}^{m-1}_{+})}^{\varepsilon_{m}} \Biggr\}^{\frac{1}{\varepsilon_{m}}}
$$
$$
\gg 2^{-n\alpha} \Biggl\{\sum\limits_{s_{m}< \frac{n}{\gamma_{m}}}V_{m}^{\varepsilon_{m}}(2^{-s_{m}})\Bigl(\prod_{j=1}^{m-1}V_{j}(2^{-n})(n+1)^{\sum\limits_{j=2}^{m-1}\frac{1}{\varepsilon_{j}}}\Bigr)^{\varepsilon_{m}} 
\Biggr\}^{\frac{1}{\varepsilon_{m}}}
$$
$$
\gg 2^{-n\alpha}\prod_{j=1}^{m}V_{j}(2^{-n})
(n+1)^{\sum\limits_{j=2}^{m}\frac{1}{\varepsilon_{j}}}.
$$
The Lemma 2 is proved.

\begin{rem} Lemma 1 for $V_{j}(t)=1$, $t\in (0, 1]$, $\theta_{j}=\theta$, $j=1,...,m$ and $\gamma_{j}=\gamma_{j}^{'}=1$ for $j = 1, \ldots, \nu $ and $1<\gamma_{j}^{'}<\gamma_{j}$ for $j=\nu + 1, ..., m$ was proved earlier in \cite[Lemma 2]{21}.
In the case $V_{j}(t)=1$, $t\in (0, 1]$,  $1\leq \theta_{j} < \infty$, $j=1,...,m$
and $\gamma_{j}=\gamma_{j}^{'}=1$ for $j=1, \ldots, \nu$ and $1<\gamma_{j}^{'}<\gamma_{j}$ for $j=\nu + 1, ..., m $ Lemma 1 is given in \cite{13}, \cite{14}.
For $V_{j}(t)=1$, $j=1,...,m$ from Lemma 1 and Lemma 2 we obtain \cite[Lemma 2 and Lemma 3]{16}.
In the case $V_{j}(t) = (1+ |\log_{2} t|)^{\lambda_{j}}$, $j = 1, ..., m $, $ t \in ( 0, 1]$  and under certain conditions on the numbers $\lambda_{j} \in \mathbb {R}$, Lemma 1 and Lemma 2 were proved in \cite{20}.
 \end{rem}

\section{The main result}\label{sec 2}

\begin{theorem}\label{th1}
Let 
 $\bar{p}=(p_{1},\ldots , p_{m})$, $\bar{q}=(q_{1},\ldots , q_{m})$,  $\bar{\tau}^{(1)}=(\tau_{1}^{(1)},\ldots , \tau_{m}^{(1)})$, $\bar{\tau}^{(2)}=(\tau_{1}^{(2)},\ldots , \tau_{m}^{(2)})$, $\bar{\gamma}=(\gamma_{1},\ldots , \gamma_{m})$, $\bar{\gamma}^{'}=(\gamma_{1}^{'},\ldots , \gamma_{m}^{'})$, $\bar{r}=(r_{1},\ldots , r_{m})$, $\bar{\theta}=(\theta_{1},\ldots , \theta_{m})$ and 
$0< \theta_{j}\le \infty$, $1<\tau_{j}^{(1)}, \tau_{j}^{(2)} <+\infty,$
$1< p_{j} <q_{j} <+\infty,$  $r_{j}> \frac{1}{p_{j}} - \frac{1}{q_{j}}$, $\gamma_{j}=\frac{r_{j} + \frac{1}{q_{j}} - \frac{1}{p_{j}}}{r_{j_{0}} + \frac{1}{q_{j_{0}}} - \frac{1}{p_{j_{0}}}}$, $1\leq \gamma_{j}^{'}\leq \gamma_{j}$,  $j=1,...,m$ and $r_{j_{0}} + \frac{1}{q_{j_{0}}} - \frac{1}{p_{j_{0}}}=\min\{r_{j} + \frac{1}{q_{j}} - \frac{1}{p_{j}} : j=1,...,m\}$, 
 $A=\{j : \frac{\gamma_{j}}{\gamma_{j}^{'}}=1, \,\, j=1,...,m\}$, $j_{1}=\min\{j\in A\}$ and  $V_{j}^{(i)}(t)=v_{j}^{(i)}(1/t)$, $t\in (0, 1]$, $j=1,\ldots, m$,  $\bar{V}^{(i)}(t)=(V_{1}^{(i)}(t),...,V_{m}^{(i)}(t))$, $i=1, 2$.

1. If 
  $1\leq\tau_{j}^{(2)}<\theta_{j}\leq +\infty$ and functions $\frac{v_{j}^{(2)}}{v_{j}^{(1)}}\in SVL[1, \infty)$, $j=1,...,m$ then
 $$
E_{n}^{(\bar{\gamma}^{'})}\left(S^{\bar r}_{\bar{p}, \bar{V}^{(1)}, \bar{\tau}^{(1)}
, \bar{\theta}}B\right)_{\bar{q}, \bar{V}^{(2)}, \bar{\tau}^{(2)}}
\asymp
2^{^{^{-n\bigl(r_{j_{0}} +\frac{1}{q_{j_{0}}}
-\frac{1}{p_{j_{0}}}\bigr)}}}\prod_{j\in A}\frac{V_{j}^{(2)}(2^{-n})}{V_{j}^{(1)}(2^{-n})}
n^{\sum\limits_{j\in A\setminus \{j_{1}\}}(\frac{1}{\tau_{j}^{(2)}}-\frac{1}{\theta_{j}})}, 
$$
for $n \in \mathbb{N}$, such that $n > n_{0}$. 

2. If   
$1\leq \theta_{j}\leq \tau_{j}^{(2)} <+\infty$, $j=1,...,m$, then the relation
 $$
E_{n}^{(\bar{\gamma}^{'})}\left(S^{\bar r}_{\bar{p}, \bar{V}^{(1)}, \bar{\tau}^{(1)}
, \bar{\theta}}B\right)_{\bar{q}, \bar{V}^{(2)}, \bar{\tau}^{(2)}}
\asymp
2^{^{^{-n\bigl(r_{j_{0}} +\frac{1}{q_{j_{0}}}
-\frac{1}{p_{j_{0}}}\bigr)}}}\prod_{j\in A}\frac{V_{j}^{(2)}(2^{-n})}{V_{j}^{(1)}(2^{-n})}
$$
holds in cases $\frac{v_{j}^{(2)}}{v_{j}^{(1)}}\in SVL[1, \infty)$, $ j = 1, ..., m $ and $A\setminus \{j_{1}\} = \emptyset$ or $\frac{v_{j}^{(2)}(t)}{v_{j}^{(1)}(t)}t^{\varepsilon}$ almost decreases for $\varepsilon> 0$ and $\frac{v_{j}^{(2)}(t)}{v_{j}^{(1)}(t)}$ is almost increasing on $ [1, \infty) $ and the set $A\setminus \{j_{1}\} \neq \emptyset$, 
for $n \in \mathbb{N}$, such that $n > n_{0}$, 
 where  $ n_ {0} $ is some positive number.
\end{theorem}

{\bf Proof.} Let 
 $1\leq\tau_{j}^{(2)}<\theta_{j}\leq +\infty$, $j=1,...,m$. 
Since 
 $\frac{v_{j}^{(2)}}{v_{j}^{(1)}} \in SVL[1, \infty )$, $j=1,\ldots, m$,  then 
  $\frac{v_{j}^{(2)}}{v_{j}^{(1)}}\in SV[1, \infty )$, $j=1,...,m$ (see [1]). Therefore, taking into account that $r_{j}> \frac{1}{p_{j}} - \frac{1}{q_{j}}$, $j = 1, \ldots, m$, we have
  $$
 \sum\limits_{s_{j}=0}^{\infty} \Bigl(2^{-s_{j}(r_{j}+ \frac{1}{q_{j}}-\frac{1}{p_{j}})}\frac{V_{j}^{(2)}(2^{-s_{j}})}{V_{j}^{(1)}(2^{-s_{j}})}\Bigr)^{\varepsilon_{j}} < \infty,
$$
where 
$\varepsilon_{j}=\tau_{j}^{(2)}\beta_{j}',$  $\frac{1}{\beta_{j}}+
\frac{1}{\beta_{j}'}=1,$  $\beta_{j}=\frac{\theta_{j}}{\tau_{j}^{(2)}}$ for 
 $j=1,\ldots, m$.
Hence, in Theorem 5 \cite{19}, setting $\varphi_{j}(t) =
t^{\frac{1}{p_{j}}}V_{j}^{(1)}(t)$ and $\psi_{j}(t) =
t^{\frac{1}{q_{j}}}V_{j}^{(2)}(t)$ for $t \in (0, 1]$ and replacing $\overline{\gamma}$ on $\overline{\gamma}^{'}$, we get
$$
E_{n}^{(\bar{\gamma}^{'})}(S_{\bar{p}, \bar{V}^{(1)}, \bar{\tau}^{(1)}, \bar\theta}^{\bar r}B)_{\bar{q}, \bar{V}^{(2)}, \bar{\tau}^{(2)}}\ll\left\|\left\{\prod_{j=1}^{m}2^{-s_{j}(r_{j}+ \frac{1}{q_{j}}-\frac{1}{p_{j}})}\frac{V_{j}^{(2)}(2^{-s_{j}})}{V_{j}^{(1)}(2^{-s_{j}})}\right\}_{\bar{s}\in Y^{m}(\bar{\gamma}^{'}, n)}
\right\|_{\bar\varepsilon}, \eqno (31) 
 $$
where 
 $\bar{\varepsilon}=(\varepsilon_{1},...,\varepsilon_{m}),$ 
$\varepsilon_{j}=\tau_{j}\beta_{j}',$  $\frac{1}{\beta_{j}}+
\frac{1}{\beta_{j}'}=1,$  $\beta_{j}=\frac{\theta_{j}}{\tau_{j}^{(2)}}.$

Now, by applying Lemma 1 with $V_{j}(t)=\frac{V_{j}^{(2)}(t)}{V_{j}^{(1)}(t)}$ and $\theta_{j}=\varepsilon_{j}$, $j = 1, \ldots, m$, $\alpha=r_{j_{0}} + \frac{1}{q_{j_{0}}} - \frac{1}{p_{j_{0}}}$ from (31), we get
  $$
E_{n}^{(\bar{\gamma}^{'})}(S_{\bar{p}, \bar{V}^{(1)}, \bar{\tau}^{(1)}, \bar\theta}^{\bar r}B)_{\bar{q}, \bar{V}^{(2)}, \bar{\tau}^{(2)}}\ll
2^{-n(r_{j_{0}} + \frac{1}{q_{j_{0}}} - \frac{1}{p_{j_{0}}})}\prod_{j\in A}\frac{V_{j}^{(2)}(2^{-n})}{V_{j}^{(1)}(2^{-n})} (n+1)^{\sum\limits_{j\in A\setminus \{j_{1}\}}(\frac{1}{\tau_{j}^{(2)}} -\frac{1}{\theta_{j}})},
$$
in the case  $1 \leq \tau_{j}^{(2)} < \theta_{j} \leq \infty$, $j=1,\ldots , m$.

Let $1 \leq \theta_{j} \leq  \tau_{j}^{(2)} < \infty$, $j = 1, \ldots, m$, then from the second statement of Theorem 5 \cite{19} and from the conditions on the functions $\frac{v_{j}^{(2)}}{v_{j}^{(1)}}$, $j=1,...,m$ follows that 
$$
E_{n}^{(\bar{\gamma}^{'})}(S_{\bar{p}, \bar{V}^{(1)}, \bar{\tau}^{(1)}, \bar\theta}^{\bar r}B)_{\bar{q}, \bar{V}^{(2)}, \bar{\tau}^{(2)}}\leq C \sup_{\bar{s} \in Y^{m}(\bar{\gamma}^{'}, n)}\prod_{j\in A}2^{-s_{j}(r_{j} + \frac{1}{q_{j}} - \frac{1}{p_{j}})}\frac{V_{j}^{(2)}(2^{-s_{j}})}{V_{j}^{(1)}(2^{-s_{j}})}
$$
$$
 \ll  2^{-n(r_{j_{0}} + \frac{1}{q_{j_{0}}} - \frac{1}{p_{j_{0}}})}\prod_{j\in A}\frac{V_{j}^{(2)}(2^{-n})}{V_{j}^{(1)}(2^{-n})}.
$$
  The upper bound is proved.
  
Let us prove the lower bound. 
Let 
 $A=\{j: \gamma_{j}^{'}=\gamma_{j}, \,\, j=1,...,m\}$, $j_{1}=\min\{j\in A\}$ and $B = \{j: \tau_{j}^{(2)} <\theta_{j},\,\, j=1,...,m\}=\{1,...,m\}$, then          
  $B^{'}=A\cap B\cup\{j_{1}\} = A\cup\{j_{1}\}=A$. We put 
   $\overline{s}^{0}=(s_{1}^{0},...,s_{m}^{0})$, where $s_{j}^{0}=s_{j}$ for $j\in A$ and $s_{j}^{0}=0$ for $j\notin A$.  

Now, consider the function
$$
f_{1, n}(\overline{x})= n^{-\sum\limits_{j\in A\setminus \{j_{1}\}}\frac{1}{\theta_{j}}}\sum\limits_{\langle\bar{s}^{0}, \bar\gamma \rangle =n}\prod_{j=1}^{m}\frac{2^{-s_{j}^{0}(r_{j}+ 1-\frac{1}{p_{j}})}}{V_{j}^{(1)}(2^{-s_{j}^{0}})}\sum_{\bar{k}\in\rho(\bar{s}^{0})}e^{i\langle\bar{k},\bar{x}\rangle},
$$
then, according to continuity, the function
 $f_{1, n}\in L_{\bar{p}, \bar{V}^{(1)}, \bar{\tau}^{(1)}}^{*}(\mathbb{T}^{m})$. 

The relation is known  (see \cite{19} ) 
$$
\Bigl
\|\sum_{\bar{k}\in\rho(\bar s)}e^{i\langle\bar{k},\bar{x}\rangle}
\Bigr\|_{\bar{p}, \bar{V}^{(1)}, \bar{\tau}}^{*}\asymp \prod_{j=1}^{m}
2^{s_{j}(1-\frac{1}{p_{j}})}V_{j}^{(1)}(2^{-s_{j}}), \eqno (32)   
 $$ 
for  $1<p_{j}, \tau_{j}<+\infty, \,\, j=1,...,m$.

Taking into account relation (32), we obtain
$$
\Bigl
\|\Bigl\{\prod_{j=1}^{m}2^{s_{j}r_{j}}\|\delta_{\bar s}(f_{1, n})\|_{\bar{p}, \bar{V}^{(1)}, \bar{\tau}^{(1)}}^{*}\Bigr\}_{\bar{s}\in \mathbb{Z}_{+}^{m}}\Bigr\|_{l_{\bar\theta}} \asymp n^{-\sum\limits_{j\in A\setminus \{j_{1}\}}\frac{1}{\theta_{j}}}\Bigl
\|\Bigl\{\chi_{\kappa(\bar{\gamma}^{'}, n)}(\bar{s}^{0})\Bigr\}_{\langle\bar{s}^{0}, \bar{\gamma}^{'} \rangle = n}\Bigr\|_{l_{\bar\theta}},   \eqno (33)  
 $$
 where  $\chi_{\kappa(\bar{\gamma}^{'}, n)}$ is the characteristic function of the set  
 $$
 \kappa(\bar{\gamma}^{'}, n)= \{\bar{s}=(s_{1},...,s_{m})\in \mathbb{Z}_{+}^{m}: \,\, \langle\bar{s}, \bar{\gamma}^{'} \rangle = n \}.
 $$
 Let 
   $\tilde{\bar{s}}=(s_{j_{1}},...,s_{j_{|A|}})$, $\tilde{\bar{\gamma}}^{'}=(\gamma_{j_{1}},...,\gamma_{j_{|A|}})$, where  $j_{i}\in A$, $i=1,...,|A|$ and $ |A|$ -- the number of elements of the set $A$, 
    then  
  $\langle\bar{s}^{0}, \bar{\gamma}^{'} \rangle = \sum\limits_{j \in A}s_{j}\gamma_{j}^{'} = \sum\limits_{i=1}^{|A|}s_{j_{i}}\gamma_{j_{i}}=\langle\tilde{\bar{s}}, \tilde{\bar{\gamma}}^{'}\rangle$.
  
  Therefore, by using Lemma 3 \cite{16}, relation (33) can be rewritten in the following form
$$
\Bigl
\|\Bigl\{\prod_{j=1}^{m}2^{s_{j}r_{j}}\|\delta_{\bar s}(f_{1, n})\|_{\bar{p}, \bar{V}^{(1)}, \bar{\tau}^{(1)}}^{*}\Bigr\}_{\bar{s}\in \mathbb{Z}_{+}^{m}}\Bigr\|_{l_{\bar\theta}} \asymp n^{-\sum\limits_{j\in A\setminus \{j_{1}\}}\frac{1}{\theta_{j}}}\Bigl
\|\Bigl\{\chi_{\kappa(\tilde{\bar{\gamma}}^{'}, |A|)}(\tilde{\bar{s}})\Bigr\}_{\langle\tilde{\bar{s}}, \tilde{\bar{\gamma}} \rangle = n}\Bigr\|_{l_{\bar\theta}} \asymp C_{1},   
 $$
 where 
  $\kappa (\tilde{\bar{\gamma}}, |A|)=\{\tilde{\bar{s}}=(s_{j_{1}},...,s_{j_{|A|}}): \,\, \langle\tilde{\bar{s}}, \tilde{\bar{\gamma}} \rangle =n \}$ and 
   $\chi_{\kappa(\tilde{\bar{\gamma}}, |A|)}(\tilde{\bar{s}})$ is the characteristic function of the set $\kappa(\tilde{\bar{\gamma}}, |A|)$.
So the function 
$F_{1, n}=C_{1}^{-1}f_{1, n} \in S_{\bar{p}, \bar{V}^{(1)}, \bar{\tau}^{(1)}, \bar\theta}^{\bar r}B$.

By the definition of the best approximation of a function and by Theorem 4 \cite{19}, we have
$$
E_{n}^{(\bar{\gamma}^{'})}(F_{1, n})_{\bar{q}, \bar{V}^{(2)}, \bar{\tau}^{(2)}} \geq \|F_{1, n}\|_{\bar{q}, \bar{V}^{(2)}, \bar{\tau}^{(2)}} 
$$
$$
\gg
\Bigl
\|\Bigl\{\prod_{j=1}^{m}2^{s_{j}(\frac{1}{\tau_{j}} - \frac{1}{q_{j}})}V_{j}^{(2)}(2^{-s_{j}})\|\delta_{\bar{s}^{0}}(f_{1, n})\|_{\bar\tau,  \bar{\tau}^{(2)}}^{*}\Bigr\}_{\langle\bar{s}^{0}, \bar{\gamma}^{'} \rangle = n}\Bigr\|_{l_{\bar{\tau}^{(2)}}}, 
$$
where
 $\bar\tau = (\tau_{1},...,\tau_{m})$ and $1< q_{j}< \tau_{j} < \infty$, $j=1,...,m$.
Hence, according to estimate (32), we obtain
$$
E_{n}^{(\bar{\gamma}^{'})}(F_{1, n})_{\bar{q}, \bar{V}^{(2)}, \bar{\tau}^{(2)}} \gg
n^{-\sum\limits_{j\in A\setminus \{j_{1}\}}\frac{1}{\theta_{j}}}
\Bigl
\|\Bigl\{2^{-\langle\bar{s}^{0}, \bar{\gamma}^{'} \rangle(r_{j_{0}}+\frac{1}{q_{j_{0}} }- \frac{1}{p_{j_{0}}})}\Bigr\}_{\langle\bar{s}^{0}, \bar{\gamma}^{'} \rangle = n}\prod_{j=1}^{m}\frac{V_{j}^{(2)}(2^{-s_{j}^{0}})}{V_{j}^{(1)}(2^{-s_{j}^{0}})}\Bigr\|_{l_{\bar{\tau}^{(2)}}} 
$$
$$
=C
n^{-\sum\limits_{j\in A\setminus \{j_{1}\}}\frac{1}{\theta_{j}}}
\Bigl
\|\Bigl\{2^{-\langle\tilde\bar{s}, \tilde\bar{\gamma} \rangle(r_{j_{0}}+\frac{1}{q_{j_{0}} }- \frac{1}{p_{j_{0}}})}\prod_{i=1}^{|A|}\frac{V_{j_{i}}^{(2)}(2^{-\tilde{s}_{j_{i}}})}{V_{j_{i}}^{(1)}(2^{-\tilde{s}_{j_{i}}})}\Bigr\}_{\tilde{\bar{s}}\in \kappa(\tilde{\bar{\gamma}}^{'}, |A|)}\Bigr\|_{l_{\bar{\tau}^{(2)}}}.
$$
Now, by using Lemma 2 for $V_{j_{i}}(t)=\frac{V_{j_{i}}^{(2)}(t)}{V_{j_{i}}^{(1)}(t)}$, $i = 1, ..., |A|$, from this, we get
$$
E_{n}^{(\bar{\gamma}^{'})}(F_{1, n})_{\bar{q}, \bar{V}^{(2)}, \bar{\tau}^{(2)}} \geq C2^{-n(r_{j_{0}}+\frac{1}{q_{j_{0}} }- \frac{1}{p_{j_{0}}})}\prod_{j\in A}
\frac{V_{j}^{(2)}(2^{-n})}{V_{j}^{(1)}(2^{-n})}(n +1)^{\sum\limits_{j\in A\setminus \{j_{1}\}}(\frac{1}{\tau_{j}^{(2)}}-\frac{1}{\theta_{j}})}
$$
in the case 
 $1\leq \tau_{j}^{(2)} < \theta_{j} \leq \infty$, $j=1,\ldots , m$.

Let 
 $1 \leq \theta_{j} \leq  \tau_{j}^{(2)} < \infty$, $j=1,\ldots , m$. 
Consider the function  
$$
f_{2, n}(\overline{x})= \prod_{j=1}^{m}\frac{2^{-s_{j}^{0}(r_{j}+ 1-\frac{1}{p_{j}})}}{V_{j}^{(1)}(2^{-s_{j}^{0}})}\sum_{\bar{k}\in\rho(\bar{s}^{0})}e^{i\langle\bar{k},\bar{x}\rangle}
$$
for  $\langle\bar{s}^{0}, \bar\gamma^{'} \rangle \geq n$.

Then, by using relation (32), we obtain that the function
$F_{2, n}=C_{1}^{-1}f_{2, n} \in S_{\bar{p}, \bar{V}^{(1)}, \bar{\tau}^{(1)}, \bar\theta}^{\bar r}B$ and 
$$
E_{n}^{(\bar{\gamma}^{'})}(F_{2, n})_{\bar{q}, \bar{V}^{(2)}, \bar{\tau}^{(2)}} \geq \|F_{2, n}\|_{\bar{q}, \bar{V}^{(2)}, \bar{\tau}^{(2)}}^{*} \gg 2^{-n(r_{j_{0}}+\frac{1}{q_{j_{0}} }- \frac{1}{p_{j_{0}}})}\prod_{j\in A}\frac{V_{j}^{(2)}(2^{-n})}{V_{j}^{(1)}(2^{-n})}.
$$
The theorem is proved.

\begin{rem} In the case $V_{j}^{(1)}(t)=V_{j}^{(2)}(t)=1$, $t\in (0, 1]$ and $p_{j}=\tau_{j}^{(1)}=p$,
$q_{j}=\tau_{j}^{(2)}=q$, $\theta_{j}=\theta$ for $j = 1,..., m$ from Theorem 1 follows the previously known results 
 by V.N. Temlyakov \cite[Theorem 2.2]{21}  and A.S. Romanyuk \cite[Theorem 2]{22}.  
For $V_{j}^{(1)}(t)=V_{j}^{(2)}(t)=1$, $t\in (0, 1]$, $j=1,...,m$ and $\gamma_{j}^{'} =\gamma_{j}=1$ for $ j = 1, ..., \nu $ and $ \gamma_{j}^{'} <\gamma_{j}$, $ j = \nu + 1, ..., m$ Theorem 1 implies \cite[Theorem 2]{13} (also see \cite[Theorem 3.5]{14}), and for $V_{j}^{(1)}(t)=V_{j}^{(2)}(t)=1$, $t\in (0, 1]$ and $\gamma_{j}^{'} \leq \gamma_{j}$ for $ j = 1, ...,m$ \cite[Theorem 1]{16}.
In the case $V_{j}^{(1)}(t)=(1+ |\log_{2} t|)^{\alpha_{j}}$ and $V_{j}^{(2)}(t)=(1+ |\log_{2} t|)^{\beta_{j}}$, $j=1,...,m$, $t\in (0, 1]$ were proved in \cite{20}.
\end{rem}

This work was supported by a grant
  Ministry of Education and Science of the Republic of Kazakhstan (Project AP 08855579).


\begin{thebibliography}{99}


\bibitem{1} Seneta E. {\it Regularly Varying Functions}. Springer-Verlag Berlin Heidelberg, 1976

\bibitem{2} 
    Bingham N., Goldie C. and Teugels J. {\it Regular Variation}, Cambridge
Univ. Press, Cambridge, 1987.

\bibitem{3}
 Edmunds~D.E., Evans~W.D. {\it Hardy operators, function spaces and embedding},  Springer-Verlag, Berlin Heidelberg - 2004. 

\bibitem{4} Stein E.M., Weiss G. {\it Introduction to Fourier analysis on
Euclidean spaces}, Princeton: Princeton Univ. Press, 1971. 

\bibitem{5} Blozinski A.P. {\it Multivariate rearragements and Banach function spaces with mixed norms}, Trans. Amer. Math. Soc.  263(1) (1981), 146-167.

\bibitem{6} Kolyada V.I. {\it Rearrangements of functions and embedding theorems}, Russian Math. Surveys,  44(5) (1989), 73-117.

\bibitem{7} Nikol'skii S. M.  {\it Approximation of functions of several variables and embedding theorems},  Moscow,  1977.

\bibitem{8} Lizorkin P.I.,  Nikol'skii S. M.  {\it Spaces of functions of mixed smoothness from the decomposition point of view}, Proc. Stekov Inst. Math.,  187 (1989), 143--161.

\bibitem{9} Amanov T.I. {\it Spaces of differentiable functions with dominant mixed derivative},   Alma-Ata,  1976.

\bibitem{10}  Tikhomirov V.M. {\it Approximation theory},
Itogy Nauki i Tekhniki :  Sovrem. Probl. Math.: Fud. Naprav. VINITI,  Moscow. 14 (1987), 103--270.
  
\bibitem{11} Dinh Dung,  Temlyakov V.N.,  Ullrich T.  {\it Hyperbolic cross approximation}, Advanced Courses in Mathematics. CRM Barcelona. Springer, Basel/Berlin. 2018. 

\bibitem{12} Temlyakov V. {\it Multivariate approximation}, Cambridge University Press. 2018. 

\bibitem{13} Akishev G. {\it Approximation of function classes in spaces with mixed norms}, Mat. Sb., 197(8) (2006), 17--40.
 
\bibitem{14}  Akishev G. {\it On approximation of function classes in Lorentz spaces with anisotropic norm}, Anal. Theory and Applic.,  29(4) (2013), 358--372.

\bibitem{15} Akishev G. {\it On Approximation Orders of Functions of Several Variables in the Lorentz Space}, Trudy Instituta Matematiki i Mekhaniki UrO RAN, 22(4) (2016) 13--28. 

\bibitem{16} Bekmaganbetov K.A. {\it On orders of approximationof the Besov class in the metric of anisotropic Lorentz spaces}, Ufim. math. jour.  1(2) (2009), 9--16.

\bibitem{17}  Bekmaganbetov K.A. {\it Orders of approximations of  Besov classes in the metric of anisotropic Lorentz spaces}, Meth. Fourier Anal. and Approx. Th., (2016), 149--158.     

\bibitem{18} Bekmaganbetov K.A., Orazgaliev E.T. {\it Bernstein-Nikol'skii inequalities and estimates of best approximation in anisotropic Lorentz spaces},  Matem. zhur., 15(2) (2015),  32--41. 

\bibitem{19} Akishev G. {\it Estimates of the order of approximation of functions of several variables in the generalized Lorentz space}, arXiv: 2105.14810v1 [mathCA] 31 may 2021. 18 p. 

\bibitem{20} Akishev G. {\it On exact estimates of the order of approximation of functions of several variables in the anisotropic Lorentz-Zygmund space}, arXiv: 2106.07188v2 [mathCA] 20 Jun 2021, 20 p.

\bibitem{21} Temlyakov V.N. {\it Approximation of functions with bounded mixed derivative}, Tr. Mat. Inst. Steklov., 178 (1986), 3--112.

\bibitem{22} Romanyuk A.S.  {\it Approximation of the Besov classes of periodic functions of several variables in the space $L_{q}$},
Ukrain . Mat. Zh., 43(10) (1991), 1297--1306.
\end{thebibliography}
 \end{document}